\documentclass[10pt]{article}
\usepackage{amsmath,amsthm,amsfonts,amssymb,amscd, amsxtra}
\usepackage[utf8]{inputenc}
\usepackage[english]{babel}
\usepackage{float}

\newtheorem{theorem}{Theorem}[section]

\newtheorem{example}{Example}[section]
\newtheorem{remark}{Remark}[theorem]

\usepackage{graphicx}

\begin{document}

\title{Newton's Method with GeoGebra}

      \author{ 
O. P. Ferreira  \thanks{IME/UFG, Campus II- Caixa Postal 131,
    CEP 74001-970 - Goi\^ania, GO, Brazil (E-mail:{\tt
      orizon@mat.ufg.br}).  The author was supported in part by
    FUNAPE/UFG,  CNPq Grant 302618/2005-8, CNPq Grant
    475647/2006-8, PRONEX--Optimization(FAPERJ/CNPq).}  
    \and
    D. A. S. Pires\thanks{IME/UFG, Campus II- Caixa Postal 131,
    CEP 74001-970 - Goi\^ania, GO, Brazil (E-mail:{\tt
      davialexandre@discente.ufg.br}). }   
         }

\date{June 16, 2021}

\maketitle

\begin{abstract}

In this work, we present a program in the computational environment, GeoGebra, that enables a graphical study of Newton's Method. Using this computational device, we will analyze Newton's Method convergence applied to various examples of real functions. Then, it will be given a guide to the construction of the program in GeoGebra

\end{abstract}
%%%%%%%%%%%%%%%%%%%%%%%
\section{Introduction}
Newton's Method is an iterative method to find numerically the solutions of nonlinear equations of the form: $f(x)=0$, with $f: D \to \mathbb{R}^n$ being differentiable and $D\subset \mathbb{R}^n$ being an open set. The idea is very simple: given a point of the domain, we will compute  the root of the linear approximation of $f$ about this point, getting a new point, for which this process will be repeated, and so on. For this process to be well defined, some hypothesis about $f$, its derivative $f'$, and the starting point $x_0$ are necessary, see \cite{DennisSchnabel1983, Ortega1970}. Moreover, in order to get results about  the convergence of the method, some additional hypotheses are necessary, see \cite{DennisSchnabel1983}. Since the method development, various results and important applications were found in several areas of pure an applied mathematics. For example, a historical  perspective  of the method's applications in optimization can be found in
 \cite{Polyak2007} and for method applications in general mathematics, see \cite{Krantz2002, Wayne1996}.  

In this paper, we present several examples of real functions and apply Newton's method to find its roots. The point is to show various situations that will help to understand the necessary hypotheses to a well  definition of the method and the formulation of convergence theorems to the sequence generated by the method. It is worth mentioning that  we use the software GeoGebra as a support tool to this study, it will be of a great help in viewing and simulating the theoretical results.

The paper presentation will be done in the following way: initially, in Section \ref{sec:gdnm}, we introduce in a formal way Newton's method. In the Section \ref{sec:exemplos}, we present several examples, highlighting analytical results through the graphic representations. In the Section \ref{sec:theos}, we present two convergence theorems and discuss the first theorem hypotheses and its relations with the previous section examples. In the Section \ref{sec:alg}, we discuss the program construction in GeoGebra Software, in ways to enable the reader to create a similar one. Finally, in Section \ref{sec:of}, we do some last remarks.
 
%%%%%%%%%%%%%%%%%%%%%%%%%%%
\section{Newton's method in $\mathbb{R}$} \label{sec:gdnm}
\begin{figure}[H] 
\centering
\includegraphics[width=0.5\textwidth]{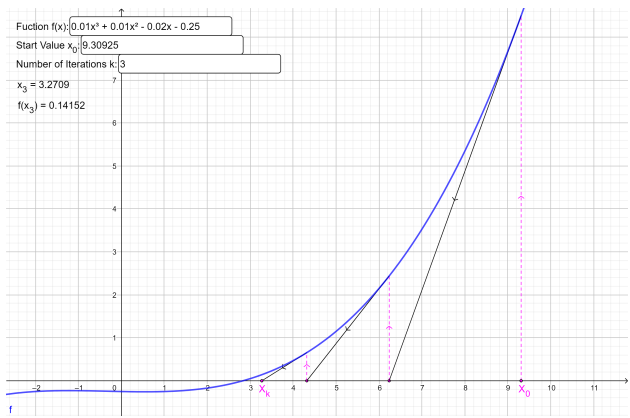}
\caption{Newton's Method} \label{fig:IMG1}
\end{figure}

In this section, we describe the general form of Newton's method applied to real functions of one variable. Let $I \subset \mathbb{R}$  be  an open interval and $f: I \to \mathbb{R}$ be a differentiable function. Take $x_0 \in I$ such that $f'(x_0) \neq 0$ and let   $ \ell(x):=f(x_0)+f'(x_0)(x-x_0)$ be the linear approximation. The point generated by Newton's method from $x_0$ is defined as the root of the equation $ \ell(x)=0$. Hence, due to   $f'(x_0) \neq 0$, we have  $x_1=x_0-f(x_0)/f'(x_0)$. If it happens that $x_i \in I$ and $f'(x_1) \neq 0$, then we repeat the process, finding a new point $x_2=x_1-f(x_1)/f'(x_1)$ and so on. This way, if for each  $n=0,1, \ldots$ we have $x_n \in I$ e $f'(x_n) \neq 0$, then we can define the Newton iteration as follows:
\begin{equation} \label{eq:nm}
x_{n+1}=x_n-\frac{f(x_n)}{f'(x_n)}, \qquad n=0,1, \ldots. 
\end{equation}
It may occur that the sequence \eqref{eq:nm} diverges or converges in several ways. In the next section, we will show various examples to understand the behavior of the sequence $(x_n)_{n\in \mathbb N}$ defined in \eqref{eq:nm}. We will now describe the general form that our examples will be shown. In Figure~\ref{fig:IMG1} we show graphically the sequence $(x_n)_{n \in \mathbb N}$ defined in \eqref{eq:nm} for the polynomial function $f(x)=0.01x^3+0.01x^2-0.02x-0.25$. 

In this figure, we can see how our examples will be generated by our program in the specific case of the polynomial function $f$. The blue curve represents the graph of the function $f$, the black lines represent the graph of the  linear approximation of  $\ell$ about the chosen point and the pink points represent each element of Newton's sequence obtained by the intersection of $\ell$ lines with the $x$-axis. As we can see, in the left superior corner of  Figure~\ref{fig:IMG1} it is shown the function, the first point, the number of iterations performed by our program, the last iteration and its image. In this case, it is important to note that only $3$ iterations of the method were represented, $x_3$ being the last point, labeled in the figure as $x_k$.
%%%%%%%%%%%%%%%%%%%%%%%%%%%
\section{Examples} \label{sec:exemplos}
In this section, we present some examples aiming to provide information to the formulation of theorems about Newton's method convergence. In all of our examples, the graphs were generated as described in Section~\ref{sec:gdnm}. First, we present an example showing that Newton's method can diverge for any choice of the initial point.
\begin{example} \label{ex:ex1}
Let  $f: \mathbb{R} \to \mathbb{R}$ be given by $ f(x)=x^{1/3}$. This function is differentiable in $ \mathbb{R}/\{0\}$ and its derivative is
$$
f'(x)=\frac{1}{3x^{2/3}},  \qquad \forall~x\ne 0.
$$
Note that  $x_*=0$ is the only root of $f$ and $\lim_{x\to 0}f'(x)=+\infty$. Some calculations shows that Newton's sequence, with initial point $x_0\neq 0$, is defined as
\begin{equation} \label{ex:den1}
x_{k+1}:=x_k-\frac{x_k^{1/3}}{1/3x_k^{2/3}}=-2x_k, \qquad k=0, 1, \ldots. 
\end{equation} 
Thus, $x_k=(-2)^kx_0$, for all  $k=0, 1, \ldots$. Therefore, we conclude that $\lim_{k\to +\infty} |x_k|=+\infty$, and then, Newton's sequence $\{x_k\}$ diverges for any choice of initial point $x_0\neq 0$. In Figure~\ref{fig:IMG2}, we graphically represent the behavior of Newton's sequence $\{x_k\}$, with initial point $x_0=0.2$, applied to the function $ f(x)=x^{1/3}$.
\begin{figure}[H] 
\centering
\includegraphics[width=0.5\textwidth]{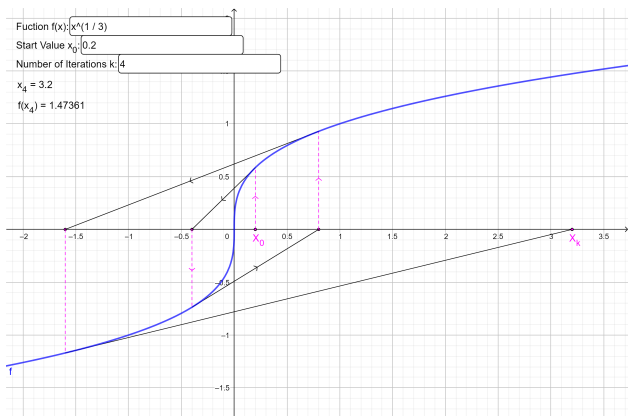}
\caption{Newton's Method applied to the function $f(x)=x^{1/3}$.} \label{fig:IMG2}
\end{figure}
\end{example}
\begin{example} \label{ex:ex2}
Let  $f: \mathbb{R} \to \mathbb{R}$ be   given by $f(x)= x^{2/3}$. Its derivative is  $f'(x)= 2/(3x^{1/3})\neq 0$, for all $x\neq 0$. Therefore, $f$ is differentiable in every point, but $x=0$. The only root of $f$  is $x=0$. Thus, given the initial point $x_0 \in \mathbb{R}/\{0\}$, the sequence generated by Newton's method is defined as
$$
x_{k+1}=x_k-\frac{x_k^{2/3}}{2/3x_k^{1/3}}=x_k-\frac{3x_k}{2}=-x_k/2, \qquad k=0, 1, \ldots. 
$$
Therefore, $x_k=(-1/2)^kx_0$, for all $k=0, 1, \ldots$. This way we conclude that the sequence $\{x_k\}$ converges, regardless of the choice of $x_0$ with linear rate  of convergence (which is slow!). In Figure~\ref{fig:IMG2.2}, we graphically represent the behavior of Newton's sequence $\{x_k\}$, with initial point $x_0=1$, applied to  $f(x)= x^{2/3}$.
\begin{figure}[H] 
\centering
\includegraphics[width=0.5\textwidth]{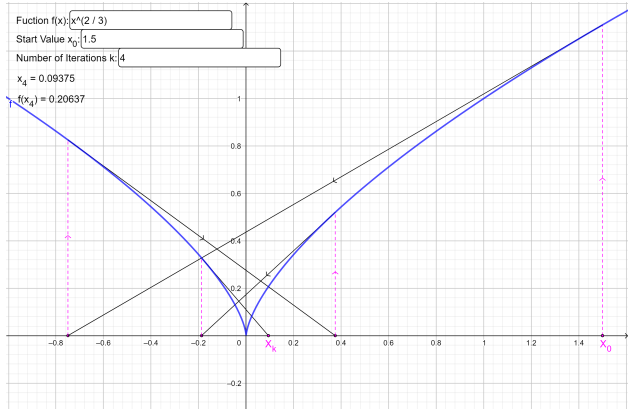}
\caption{Newton's Method applied to the function $f(x)= x^{2/3}$.} \label{fig:IMG2.2}
\end{figure}
\end{example}

\begin{example} \label{ex:ex3}
Let  $f: \mathbb{R} \to \mathbb{R}$ be  the function defined by $f(x)=x^3$, whose derivative is $f'(x)=3x^2$, which exists and it is  not $0$ whatever be the $x \neq 0$. $f$ has only one root, $x=0$, and given an initial point $x_0 \neq 0$, we get Newton's sequence:
$$
x_{k+1}=x_k-\frac{x_k^3}{3x_k^2}=x_k-\frac{x_k}{3}=\frac{2x_k}{3}, \qquad k=0, 1, \ldots, 
$$
which give us $x_k=(2/3)^kx_0$, for every $k=0,1, \ldots$. Such sequence converges, whatever be the chosen $x_0$. One can see the graphical representation of the method to the function $f$, with initial point $x_0=1$, in Figure~\ref{fig:IMG3}.

\begin{figure}[H] 
\centering
\includegraphics[width=0.5\textwidth]{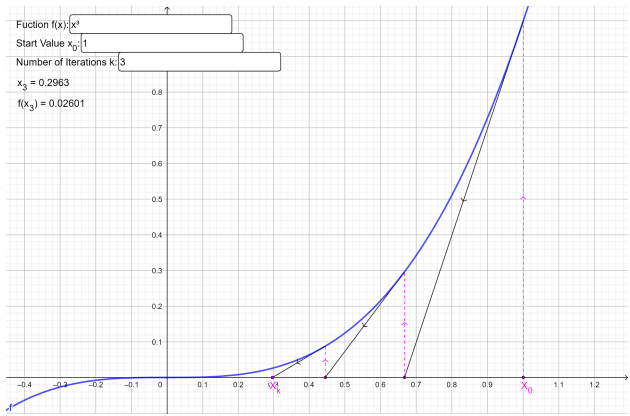}
\caption{Newton's method applied to the function $f(x)= x^3$.} \label{fig:IMG3}
\end{figure}
\end{example}

\begin{example} \label{ex:ex4}
Let  $f: \mathbb{R} \to \mathbb{R}$ be  such that $f(x)= x/\sqrt{1+x^2}$, which has  $x=0$  as its unique root  and  derivative $f'(x)=1/(1+x^2)^{3/2}\neq 0$  for all $x\in \mathbb{R}$. Thus, given a point $x_0 \in \mathbb{R}$, Newton's sequence will be:
$$
x_{k+1}=x_k-\frac{x_k/\sqrt{1+x_k^2}}{1/(1+x_k^2)^{3/2}}=x_k-x_k(1+x_k^2)=-x_k^3
$$
This way, we have $x_k=(-1)^kx_0^{3^k}$, which defines a sequence in three cases:
\begin{description}
\item [(i)] If $|x_0|<1$, then the sequence converges to $0$. In Figure \ref{fig:IMG5}, one can see the method with $x_0=0.9$ and $3$ iterations.
\begin{figure}[H]
\centering
\includegraphics[width=0.5\textwidth]{IXP51.pdf}
\caption{Newton's Method applied to function $f(x)= x/\sqrt{1+x^2}$, with $|x_0|<1$} \label{fig:IMG5}
\end{figure}
\item [(ii)] If $|x_0|=1$, then the sequence oscillates between $-1$ and $1$.  In Figure \ref{fig:IMG6}, one can see the method with $x_0=1$ and $3$ iterations.
\begin{figure}[H]
\centering
\includegraphics[width=0.5\textwidth]{IXP52.pdf}
\caption{Newton's Method applied to function $f(x)= x/\sqrt{1+x^2}$, with $|x_0|=1$} \label{fig:IMG6}
\end{figure}

\item [(iii)] If $|x_0|>1$, then the sequence diverges. In Figure \ref{fig:IMG7}, one can see the method $x_0=1.02$ and $3$ iterations.
\end{description}

\begin{figure}[H]
\centering
\includegraphics[width=0.5\textwidth]{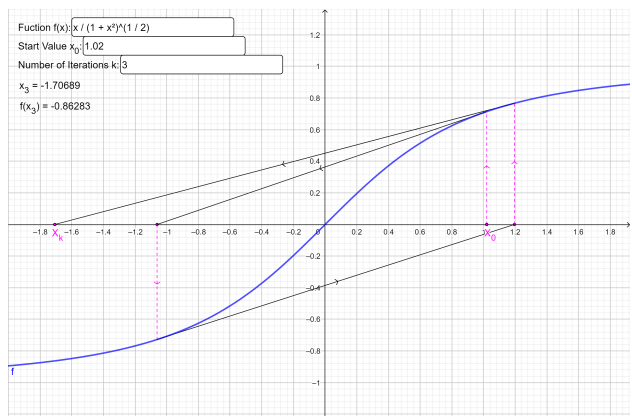}
\caption{Newton's Method applied to function $f(x)= x/\sqrt{1+x^2}$, with $|x_0|>1$} \label{fig:IMG7}
\end{figure}
\end{example}

\begin{example} \label{ex:ex5}
Consider $f: (0,+\infty) \to \mathbb{R}$, defined by $f(x)=1-1/x$, whose derivative is $f'(x)=1/x^2$. Note that $f$ has  a unique  root  $x=1$. Given an initial point $x_0$, the Newton's sequence is  as follows:
$$
x_{k+1}=x_k-\frac{1-1/x_k}{1/x_k^2}=x_k-x_k^2+x_k=x_k(2-x_k)
$$
Performing some algebraic manipulations on the last equation, we can show that $1-x_{k+1}=(1-x_k)^2,$ for all $k=0,1, \ldots.$ Using induction argument, we can show that:
\begin{equation} \label{eq:exe}
x_{k}=1-(1-x_0)^{2^k}, \qquad k=0,1, \ldots.
\end{equation}
Note that there are two cases to consider:
\begin{description}
\item [(i)] If it is  the case of $x_0 \geq 2$, then $x_1 \leq 0$, which is not a point of the domain of $f$. Thus, Newton's sequence is not well defined to this choice of  $x_0$, see Figure~\ref{fig:IMG8}.

\begin{figure}[H]
\centering
\includegraphics[width=0.5\textwidth]{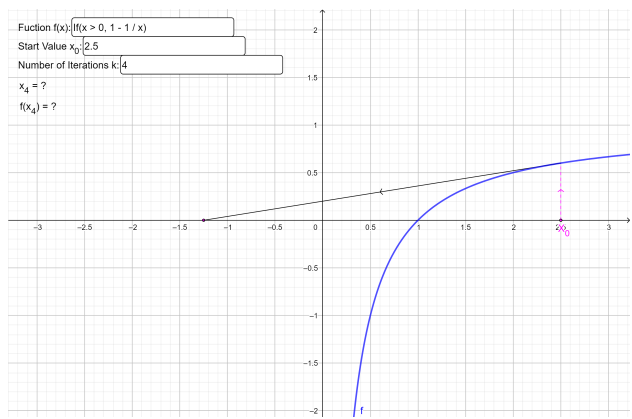}
\caption{Newton's Method applied to the function $f(x)= 1-1/x; x \in (0,+\infty)$, with $x_0 \geq 2$} \label{fig:IMG8}
\end{figure}

\item [(ii)] If  $0<x_0<2$, take the limit on \eqref{eq:exe} to conclude that $x_k$ converges to $1$, which is a root of $f$, see Figure~\ref{fig:IMG9}.
\end{description}

\begin{figure}[H]
\centering
\includegraphics[width=0.5\textwidth]{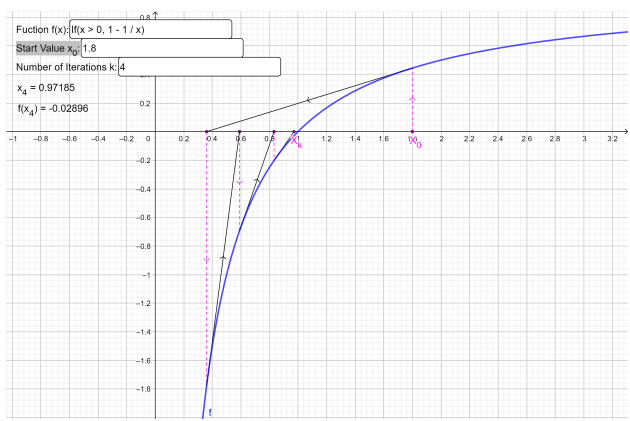}
\caption{Newton's Method applied to the function $f(x)= 1-1/x; x \in (0,+\infty)$, with $0<x_0<2$} \label{fig:IMG9}
\end{figure}
\end{example}

\begin{example} \label{ex:ex6}
Define  $f: \mathbb{R} \to \mathbb{R}$ as $f(x)=x(x-1)(x+1)=x^3-x$. It  is immediate that the roots of $f$ are $-1,0$  and $1$.  The derivative of $f$ is $f'(x)=3x^2-1$, which is $0$ in the points $-1/\sqrt{3} $ and $1/\sqrt{3}$. Newton's sequence is given by:
\begin{equation} \label{eq:niqp}
x_{k+1}=x_k-\frac{x_k^3-x_k}{3x_k^2-1}=\frac{3x_k^3-x_k-x_k^3+x_k}{3x_k^2-1}=\frac{2x_k^3}{3x_k^2-1},  \qquad k=0, 1, \ldots.
\end{equation}
For this function, there are three interesting cases to be considered:
\begin{description}
\item[(i)] If $|x_0|=1/\sqrt{5}$, the sequence oscillates between $1/\sqrt{5}$ and $-1/\sqrt{5}$. Indeed, replacing $x_0=1/\sqrt{5}$ in \eqref{eq:niqp} gives us:
$$
x_{1}=\frac{2x_0^3}{3x_0^2-1}=\frac{2(5^{-3/2})}{3(5^{-1})-1}=-\frac{1}{\sqrt{5}}, \qquad x_{2}=\frac{2x_1^3}{3x_1^2-1}=\frac{2(-5^{-3/2})}{3(5^{-1})-1}=\frac{1}{\sqrt{5}}.
$$
This is shown in Figure~\ref{fig:IMG10}.
\begin{figure}[H]
\centering
\includegraphics[width=0.5\textwidth]{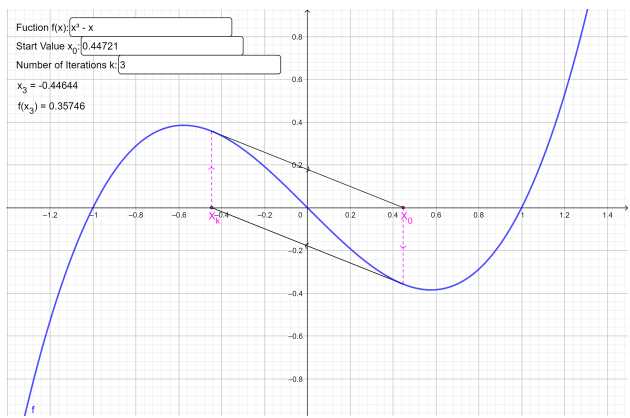}
\caption{Newton's Method applied to the function $f(x)= x^3-x$, with $|x_0|=1/\sqrt{5}$} \label{fig:IMG10}
\end{figure}
\item[(ii)] If $|x_0|=1/2$, Newton's sequence converges in one iteration. Indeed, replacing $x_0=1/2$ in \eqref{eq:niqp} gives us:
$$
x_{1}=\frac{2x_0^3}{3x_0^2-1}=\frac{2(2^{-3})}{3(2^{-2})-1}=-1.
$$
This is shown in Figure~\ref{fig:IMG11}.
\begin{figure}[H]
\centering
\includegraphics[width=0.5\textwidth]{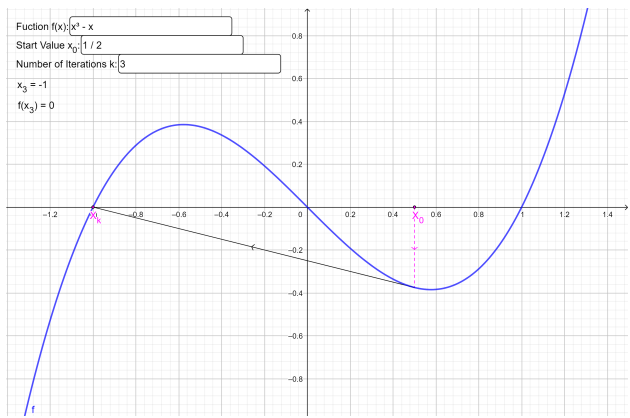}
\caption{Newton's Method applied to the function $f(x)= x^3-x$, with $|x_0|=1/2$} \label{fig:IMG11}
\end{figure}
\item[(iii)] If $x_0$ is such that $|x_1| =1/\sqrt{3}$, then $x_1$ is a root of the derivative and Newton's sequence is not well defined. Note that   such values of $x_0$, namely $1/\sqrt{3}$ and $-1/\sqrt{3}$, are real roots for the polynomials $p(x)=2\sqrt{3}x^3-3x^2+1$ and $q(x)=2\sqrt{3}x^3+3x^2-1$ respectively, which are the points $-0.4656$ and $0,4656$, rounded to 4 decimal places. We exemplify this fact in Figure~\ref{fig:IMG12} with $x_0=0,4656...$.
\end{description}

\begin{figure}[H]
\centering
\includegraphics[width=0.5\textwidth]{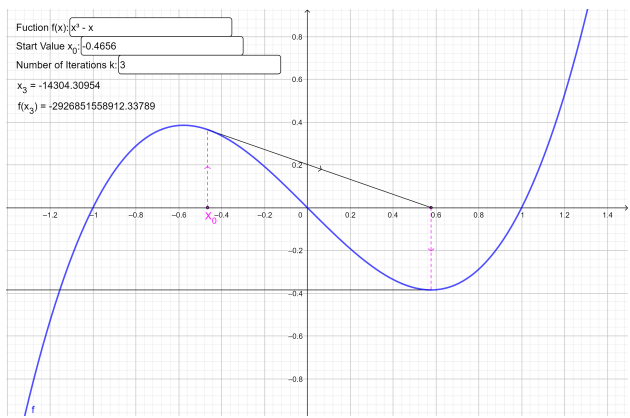}
\caption{Newton's Method applied to the function $f(x)= x^3-x$, with $|x_0|=0.4656$} \label{fig:IMG12}
\end{figure}
Note that, because it is an approximation, Figure \ref{fig:IMG12} still generates a Newton's sequence, with an extremely far $x_1$ element. Thus, the more refined is the approximation, the more distant will be the point $x_1$, in such a way that, if we take the exact root of the polynomial $2\sqrt{3}x^3-3x^2+1$, then the tangent line to $f$ in the point $(x_1,f(x_1))$ will be parallel to the $x$-axis.

\end{example}

\begin{example} \label{ex:ex7}
Let $f: \mathbb{R}/\{0\} \to \mathbb{R}$ be defined by $f(x)=|x|^x+e^x+\ln|x|+x^{(1/3)}$. Note that it is not easy to obtain the roots of $f$ without the aid of computers. The same occurs with its graph. Thus, let us  not focus so much on the algebraism behind this function. What we want to point out  is illustrated in Figure \ref{fig:IMG13}: even though $x_0=-0.65$ is a point between two roots of $f$, ($-0.19896$ and $-1.55034$, rounded to 5 decimal places), the method converges to another root, far away from $x_0$, $-6.37706$. We can take this anomaly to the extreme by taking $x_0=-0.6$ and applying the method with 6 iterations, checking that the method converges to an even more distant root, $-93.35446$. However, we will not  show a figure of this last case, because it is too complicated given the distances of the graphical elements.

\begin{figure}[H]
\centering
\includegraphics[width=0.5\textwidth]{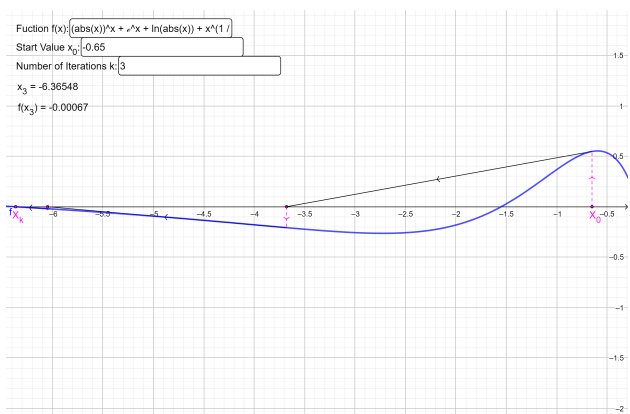}
\caption{Newton's Method applied to the function $f(x)= |x|^x+e^x+\ln|x|+x^{(1/3)}$, with $x_0=-0.65$} \label{fig:IMG13}
\end{figure}
\end{example}

\section{Convergence theorems} \label{sec:theos}
In this section, we will formulate two theorems that explain the convergence of the illustrated method on the examples of the previous section.
\begin{theorem}[Basic Theorem on Convergence] \label{th:tb}
Let $I$ be an open interval, the real valued function $f: I \to \mathbb{R}$ be  continuously differentiable and $x_* \in I$. If it is the case of $f(x_*)=0$ and $f'(x_*) \neq 0$ then, there is a $\delta>0$ such that for all $x_0 \in (x_*-\delta, x_*+\delta) \subset I$, Newton's sequence given by:
$$
x_{k+1}=x_k-\frac{f(x_k)}{f'(x_k)} \qquad k=0,1,2, \ldots
$$
is well defined, completely contained in $(x_*-\delta, x_*+\delta)$ and converges to point $x_*$. Furthermore, the rate of convergence is superlinear, as it follows:
$$
\lim_{k \to \infty}\frac{|x_{k+1}-x_*|}{|x_k-x_*|}=0
$$
\end{theorem}
Next, we will discuss some aspects of Theorem~\ref{th:tb} above considering the examples shown in Section~\ref{sec:exemplos}, this will help us to understand its hypotheses and thesis.
\begin{remark}
The following items refer to the examples covered in Section~\ref{sec:exemplos}.
\begin{description}
\item[(1)] Clearly, Theorem~\ref{th:tb} above does not apply to Example~\ref{ex:ex1}. Indeed, the function is not differentiable on its root. However, one can see that Theorem~\ref{th:tb} hypotheses are not necessary to Newton's sequence to be well defined. In this case, the sequence diverges, as it can be seen in \eqref{ex:den1} and in Figure~\ref{fig:IMG2} too.
\item[(2)]  Theorem~\ref{th:tb} does not apply to Example~\ref{ex:ex2} either, this is due to the fact that the function is not differentiable on its root. In this case, one can see that besides Newton's sequence being well defined for any initial point other than $0$, it converges to the root. However, the rate of convergence of the sequence is only linear, see Figure~\ref{fig:IMG2.2}. 
\item[(3)] The function shown in Example~\ref{ex:ex3} is differentiable, however Theorem~\ref{th:tb} does not apply due to the derivative of the function being $0$ in the root. In this case, one can see that Newton's sequence is well-defined to any initial point other than $0$, and it converges linearly to the root, see Figure~\ref{fig:IMG3}.
\item[(4)] In Example~\ref{ex:ex4}, Theorem~\ref{th:tb} applies perfectly. In this case, $\delta=1$ and $f'(0)=-1$. Note that for any initial point out of the interval $(-1, 1)$, the sequence diverges. It is worth noting that the rate of convergence is cubic if the initial point is in the interval $(-1, 1)$, this show us that the rate of convergence can be faster than superlinear rate, see Figure~\ref{fig:IMG5}. 
\item[(5)] It is visible that Theorem~\ref{th:tb} does apply to the function of Example~\ref{ex:ex5}. For this case, $\delta=1$ e $f'(1)=1$. Note that, for any initial point $x_0 \geq 2$, Newton's sequence is not well defined due to the fact that the iterate does not belong to the domain of $f$, that is $x_1\leq 0$. It is worth saying that the rate of convergence is quadratic for all initial points belonging to the interval $(0,2)$, see Figure~\ref{fig:IMG9}.
\item[(6)] One can note that, for Example~\ref{ex:ex6}, Theorem~\ref{th:tb} is applicable to the three roots of the function $f$, since $f'(-1)=2$, $f'(0)=1$ and $f'(1)=2$. We can deduce that, for the root $x_*=0$, $\delta=1/\sqrt{5}$ and, by Theorem~\ref{th:tb}, there is a maximum interval around each one of the other two roots in which not only the convergence is guaranteed, but also it has a superlinear rate of convergence, at least.
\item[(7)] It is not algebraically clear that the Theorem~\ref{th:tb} holds for Example~\ref{ex:ex7}. However, one can see through Figure~\ref{fig:IMG13} that the method does apply to the 3 roots shown, since none of them occur in a point whose tangent's inclination is $0$. The same can be said about the other two roots not shown. The most important is to note that, even though the point $x_0=-0.65$ is closer to the root $x_*=-0.19896$, the method converges to other root, $x_*=-6.37706$, which indicates us that the convergence interval can be very small. Furthermore, the interval of convergence around the root given by Theorem~\ref{th:tb} is not necessarily the biggest set of initial points for which the method converges to this specific root. Note that this is an example of difficult analysis and that the program does all the work of applying the method graphically. In the next section, we will explain the created program's algorithm, so that the interested reader may reproduce it and build its own examples. Furthermore, one can have access to the built program in GeoGebra's site: {\rm https://www.geogebra.org/m/j9hx3abd}.
\end{description}
\end{remark}
The next theorem show us that we can add an hypothesis, obtaining the greatest possible convergence interval around a root by doing so, \cite{Wang1999}.
\begin{theorem}[Convergence Theorem under Lipschitz's Condition]
Let  $I$ be  an open interval, a function $f: I \to \mathbb{R}$ be  continuously differentiable and $x_* \in I$. Now, suppose that $f(x_*)=0$, $f'(x_*) \neq 0$ and that there is $K>0$ such that:
$$
\frac{|f'(x)-f'(y)|}{|f'(x_*)|} \leq K|x-y|, \qquad \forall x,y \in I
$$
Let be $\kappa:=\sup \{t>0: (x_*-t,x_*+t) \subset I\}$ and $r:=\min \{ \kappa, 2/(3K) \}$. Thus, Newton's sequence given by:
$$
x_{k+1}=x_k-\frac{f(x_k)}{f'(x_k)} \qquad k=0,1,2, \ldots
$$
is well defined, completely contained in $(x_*-\delta, x_*+\delta)$, converges to the point $x_*$, and the following inequality holds:
$$
|x_*-x_{k+1}| \leq \frac{K}{2(1-K|x_0-x_*|)} |x_k-x_*|^2 \qquad k=0,1,2, \ldots
$$
Furthermore, $x_*$ is the only root of $f$ in the open interval $(x_*-2/K, x_*+2/K)$ and if $2/(3K) < \kappa$, then $(x_*-2/(3K),x_*+2/(3K))$ is the greatest possible convergence interval around $x_*$.
\end{theorem}
%%%%%%%%%%%%%%%%%%%%%%%%
\section{Programming on GeoGebra} \label{sec:alg}
In this section, it will be explained how the used program was developed to build the figures above, in a way that the reader may build a similar one, and through it, analyze its own functions with Newton's method, having control over the initial point, the number of iterations and the aesthetic elements. The program was developed with GeoGebra Classic, and it can be accessed through the link: {\rm https://www.geogebra.org/m/j9hx3abd}.
\subsection{Defining the functions}
\begin{figure}[H]
\centering
\includegraphics[width=0.5\textwidth]{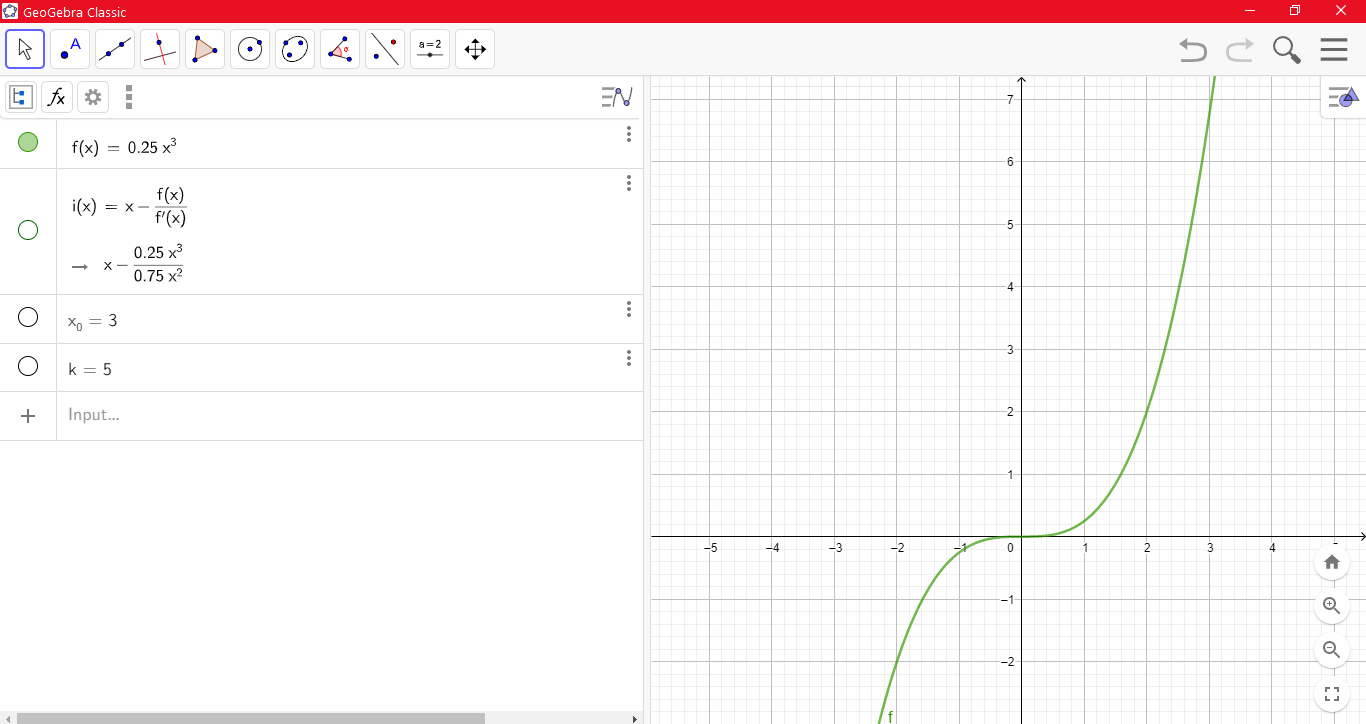}
\caption{Defining Functions and Constants} \label{fig:IMG14}
\end{figure}

The first step is to define the algebraic functions. Define the function to be analyzed, $f(x)$ and the iteration function as $i(x)=x-f(x)/f'(x)$ (do not write $f$ function explicitly on $i$'s definition, so that $i$ may keep linked to $f$ in such a way that if $f$ function changes, $i$ will also change automatically). Up next, define $x_0$ and $k$, the initial point and the number of iterations, respectively, and specify in  its configurations that the object $k$ must be a natural number (set its increment to $1$ and its minimum value to $1$). Lastly, erase the unnecessary elements of the graphic interface (such as $i$'s graphic and the sliding controls of $x_0$ and $k$ in the algebraic window), making them invisible. At the end of this step, algebra and graphic windows should be looking like  Figure~\ref{fig:IMG14}.

\subsection{Defining the lists}
\begin{figure}[H]
\centering
\includegraphics[width=0.5\textwidth]{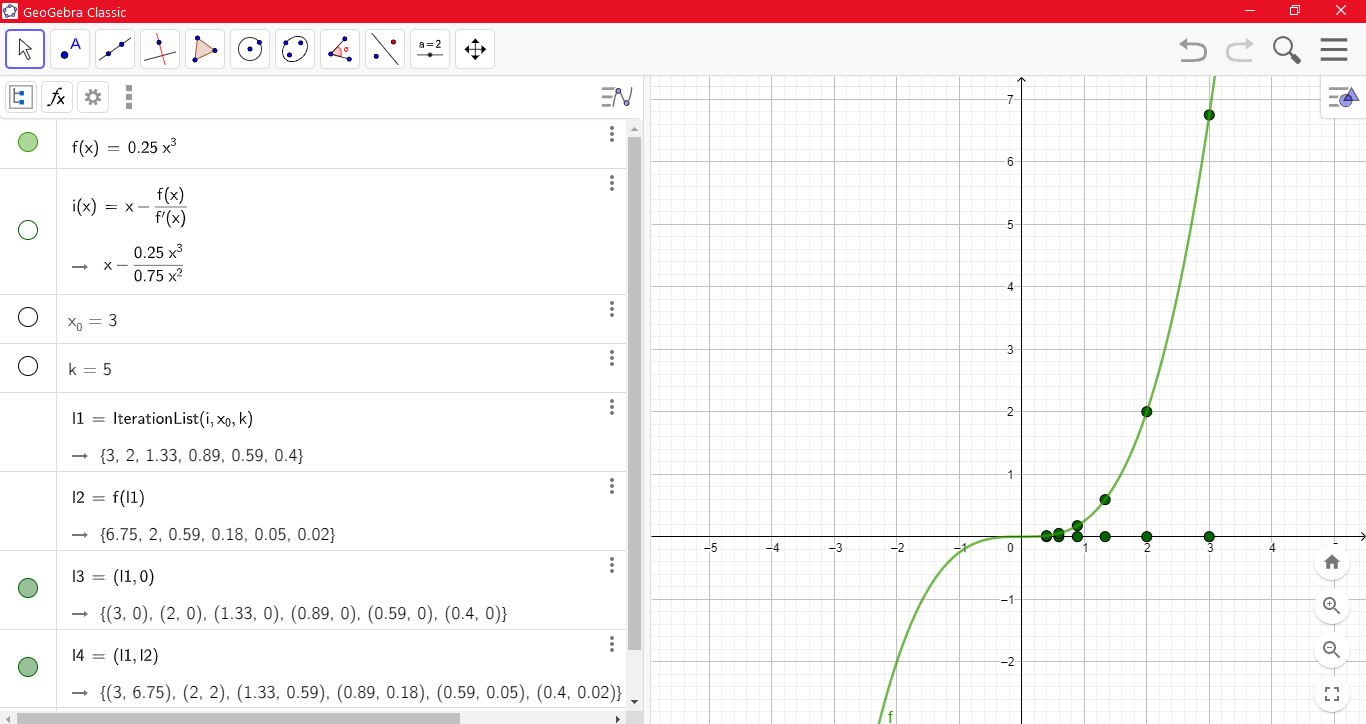}
\caption{Defining the Lists} \label{fig:IMG15}
\end{figure}

Now, we will define several lists (both of numbers and points), that will be necessary to the method's representation. Using the command  ``IterationList'', define the list $l1$ in algebra window as it follows: $l1=$ IterationList$(i, x_0, k)$ This command generates a list that applies the $i$ function to the initial value, $x_0$, $k$ times, and it begins with $x_0$, showing the successive applications, being a list of $k+1$ elements. Up next, define the $l2$ list as it follows: $l2=f(l1)$, that is, a list of images of the $l1$ list by the $f$ function. Lastly, define the point lists $l3$ and $l4$ by: $l3=(l1,0)$ e $l4=(l1,l2)$; thus, $l3$ shows the points of Newton's sequence on the x-axis and $l4$ shows the images of these points on the graph of $f$ function. At the end of this step, algebra and graphic windows should be looking like Figure~\ref{fig:IMG15}.

\subsection{Defining the line segments}
\begin{figure}[H]
\centering
\includegraphics[width=0.5\textwidth]{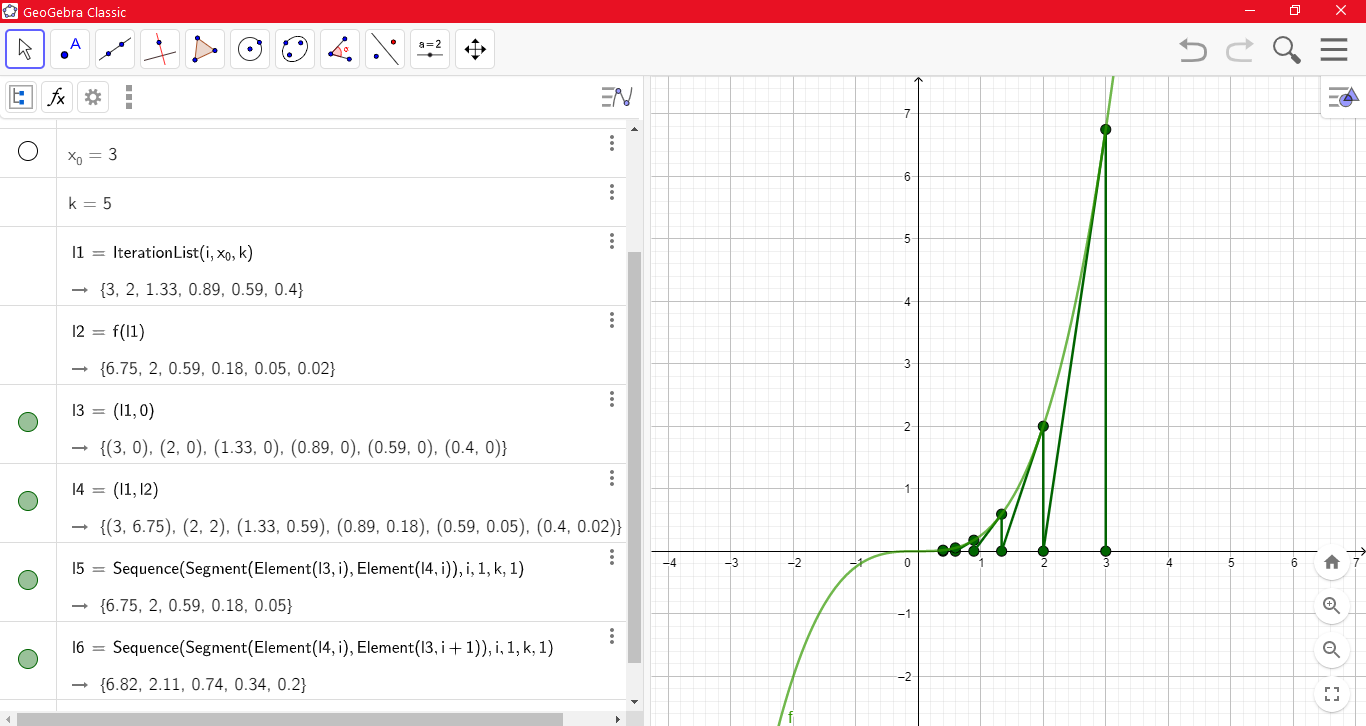}
\caption{Defining the Segments} \label{fig:IMG16}
\end{figure}

The last elements to be defined are the line segments that shows us exactly how the method works graphically (the algebraic part, Newton's sequence itself, that was defined in the first step is not enough to see how the points of such sequence were built). Here, we will need the ``Sequence'' command to define segment lists ( using the ``Segment'' command and taking lists instead of points does not work). The used command, ``Sequence'', runs other command several times and takes 5 parameters: ``Expression'': the command to be performed,  ``Variable'': a control variable, ``Start Value'': the first assignment to the control variable, ``End Value'': the maximum value to the control variable and ``Increment'': the value to be incremented to the control variable for each running of the specified command in ``Expression''. The performed command will be``Segment'', that draws a line segment and takes two parameters: two points on the plane, that will be the ends of the line segment. However, note that we have not defined any isolated point until now, only two point lists $l3$ and $l4$. To call the points of these lists, we need another command: ``Element'', that specifies an element of a list and takes two parameters: a list and a natural number, that represents the desired element of the specified list.
Thereby, define the segment list, $l5$ as being: $l5=$ Sequence(Segment(Element($l3,i$), Element($l4,i$)), $i$, $1$, $k$, $1$).Thus, $l5$ will be a list of $k$ elements, line segments between the points of Newton's sequence and its images by the $f$ function, except for the last point (remember that Newton's sequence element list has $k+1$ points, starting with $x_0$ and ending with $x_k$).
Now, define the segment list, $l6$ with the same command, just changing the points: $l6=$ Sequence(Segment(Element($l4,i$), Element($l3,i+1$)), $i$, $1$, $k$, $1$). $l6$ is a list of segments between the image of each element of Newton's sequence and the next element on the x-axis. These segments are, by definition, tangent to the graph of $f$ function. In algebra window, one can see the lists $l5$ and $l6$ as lists of $k$ real numbers; these numbers are the lengths of the segments that represent each element of the sequence and have no importance at all. At the end of these steps, the environment should be looking like Figure~\ref{fig:IMG16}.

\subsection{Adjusting the aesthetic elements}
\begin{figure}[H]
\centering
\includegraphics[width=0.5\textwidth]{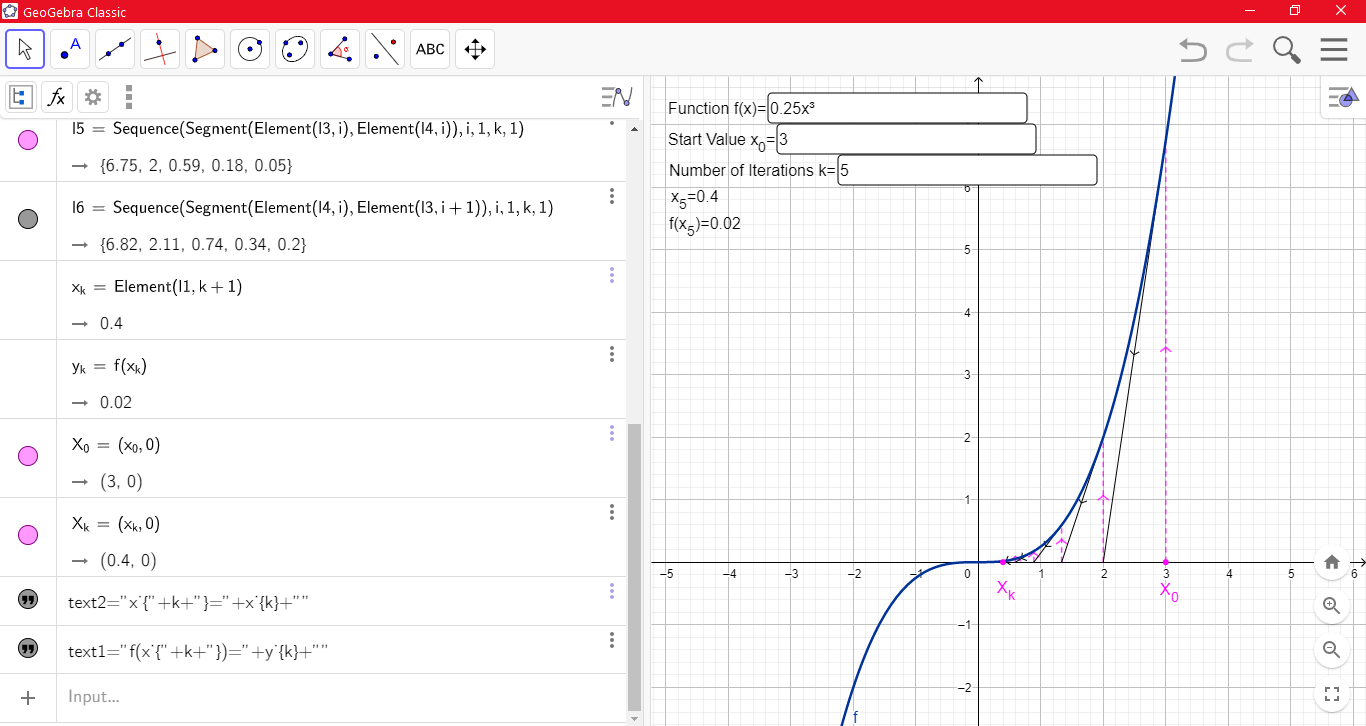}
\caption{Finishing the Program} \label{fig:IMG17}
\end{figure}

The whole structure of the program is ready but, one can see that, the way it presents itself like in Figure~\ref{fig:IMG16}, the program is not only inelegant but low functional too, since one cannot move the initial point, dragging it continuously, only changing it value manually. Moving the initial point over the x-axis gives us a dynamical idea of the method's convergence (or divergence). All of the following steps are optional, but we advise the reader to follow the steps below in order to improve its experience with the program. To modify any element's aesthetics, just select it on algebra window, open the element's configurations and go to the ``Style'' tab).
Increase the thickness of $f$'s function line, to highlight it among the other segments; reduce the size or hide the points of the lists $l3$ and $l4$; reduce the segments' thickness of $l5$ and $l6$ lists and modify its styles, so that the vertical and oblique lines look different from each other and add a arrow to the segments, this helps us to follow the path created by the method from the initial point to the final point. Define a point $X_0=$ ($x_0,0$): being linked to the $x_0$ constant, and separated from the point lists, this point can be moved, shifting the position of all the other elements in the figure. Define the constants, $x_k$ and $y_k$, by $x_k=$ Element($l1, k+1$) and $y_k=f(x_k)$. Define the point $X_k=(x_k,0)$ too (this point is only being defined to generate a label to the last point in Newton's sequence). Change the label's positions of the points $X_0$ and $X_k$ as you wish.
Make some text boxes on the graphical window, which allow the user to modify the function, $x_0$ and $k$, besides of  giving a more explanatory subtitle than the algebra window. Since we are working with Optimization and Numerical Calculus, it is interesting to modify GeoGebra's configurations so that it shows more decimal places of the represented numbers. A possible final looking is the one that follows in Figure~\ref{fig:IMG17}.

%%%%%%%%%%%%%%%%%%%%%%%%%%%%%%%%
\section{Final considerations} \label{sec:of}
The program described on Section~\ref{sec:alg} shows us that we can move the initial point in a dynamical way, observing the behavior of the sequence generated by Newton's method. The developed program aims to assist the study of Newton's method only to real valued functions. It would be interesting to build a program that aids the study the  behavior  of the Newton method to solve nonlinear equations higher dimensions.

%\bibliographystyle{abbrv}
%\bibliography{NewtGeo}

\end{document}